\documentclass[a4paper,12pt]{article}
\usepackage{amsmath}
\usepackage{amsthm}
\usepackage{amssymb}
\usepackage{enumerate}
\usepackage{url}
\usepackage[utf8]{inputenc}
\usepackage{comment}

\title{On two-dimensional extensions of Bougerol's identity in law}

\author{Yuu Hariya\thanks{{\it Corresponding author. E-mail:~hariya@tohoku.ac.jp}} \thanks{Mathematical Institute, Tohoku University, Aoba-ku, Sendai 980-8578, Japan.} \and Yohei Matsumura\thanks{Ministry of Internal Affairs and Communications, Chiyoda-ku, Tokyo 100-8926, Japan.}}

\date{\empty}
\numberwithin{equation}{section}

\theoremstyle{plain}

\newtheorem{thm}{Theorem}[section]

\newtheorem{lem}{Lemma}[section]

\theoremstyle{definition}

\theoremstyle{remark}
\newtheorem{rem}{Remark}[section]

\begin{document}

\newcommand\ND{\newcommand}
\newcommand\RD{\renewcommand}

\ND\N{\mathbb{N}}
\ND\R{\mathbb{R}}
\ND\Q{\mathbb{Q}}
\ND\C{\mathbb{C}}

\ND\F{\mathcal{F}}

\ND\kp{\kappa}

\ND\ind{\boldsymbol{1}}

\ND\al{\alpha }
\ND\la{\lambda }
\ND\La{\Lambda }
\ND\ve{\varepsilon}
\ND\Om{\Omega}

\ND\ga{\gamma}
\ND\be{\beta}

\ND\lref[1]{Lemma~\ref{#1}}
\ND\tref[1]{Theorem~\ref{#1}}
\ND\pref[1]{Proposition~\ref{#1}}
\ND\sref[1]{Section~\ref{#1}}
\ND\ssref[1]{Subsection~\ref{#1}}
\ND\aref[1]{Appendix~\ref{#1}}
\ND\rref[1]{Remark~\ref{#1}} 
\ND\cref[1]{Corollary~\ref{#1}}
\ND\eref[1]{Example~\ref{#1}}
\ND\fref[1]{Fig.\ {#1} }
\ND\lsref[1]{Lemmas~\ref{#1}}
\ND\tsref[1]{Theorems~\ref{#1}}
\ND\dref[1]{Definition~\ref{#1}}
\ND\psref[1]{Propositions~\ref{#1}}
\ND\rsref[1]{Remarks~\ref{#1}}
\ND\sssref[1]{Subsections~\ref{#1}}

\ND\pr{\mathbb{P}}
\ND\ex{\mathbb{E}}
\ND\br{W}

\ND\prb[2]{P^{(#1)}_{#2}}
\ND\exb[2]{E^{(#1)}_{#2}}

\ND\eb[1]{e^{B_{#1}}}
\ND\ebm[1]{e^{-B_{#1}}}
\ND\hbe{\Hat{\beta}}
\ND\hB{\Hat{B}}
\ND\argsh{\mathrm{Argsh}\,}
\ND\zmu{z_{\mu}}
\ND\GIG[3]{\mathrm{GIG}(#1;#2,#3)}
\ND\gig[3]{I^{(#1)}_{#2,#3}}
\ND\argch{\mathrm{Argch}\,}

\ND\vp{\varphi}
\ND\eqd{\stackrel{(d)}{=}}
\ND\db[1]{B^{(#1)}}
\ND\da[1]{A^{(#1)}}
\ND\dz[1]{Z^{(#1)}}
\ND\Z{\mathcal{Z}}

\ND\anu{\alpha ^{(\nu )}}

\ND\Ga{\Gamma}
\ND\calE{\mathcal{E}}
\ND\calD{\mathcal{D}}

\ND\f{F}
\ND\g{G}

\ND\tr{\mathbb{T}}
\ND\T{T}

\ND{\rmid}[1]{\mathrel{}\middle#1\mathrel{}}

\def\thefootnote{{}}

\maketitle 
\begin{abstract}
Let $B=\{ B_{t}\} _{t\ge 0}$ be a one-dimensional standard 
Brownian motion and denote by $A_{t},\,t\ge 0$, the quadratic 
variation of $e^{B_{t}},\,t\ge 0$. 
The celebrated Bougerol's identity in law (1983) asserts that, if 
$\beta =\{ \beta _{t}\} _{t\ge 0}$ is another Brownian motion 
independent of $B$, then $\beta _{A_{t}}$ has the same law as 
$\sinh B_{t}$ for every fixed $t>0$. Bertoin, Dufresne and Yor 
(2013) obtained a two-dimensional extension of the identity 
involving as the second coordinates the local times of $B$ and 
$\beta $ at level zero. In this paper, we present a generalization 
of their extension in a situation that the levels of those 
local times are not restricted to zero. Our argument 
provides a short elementary proof of the original extension and 
sheds new light on that subtle identity.
\end{abstract}

\section{Introduction}\label{;intro}

This paper concerns a two-dimensional extension of Bougerol's 
identity in law obtained by Bertoin, Dufresne and Yor \cite{bdy}.
Let $B=\{ B_{t}\} _{t\ge 0}$ be a one-dimensional standard 
Brownian motion and set $A_{t},\,t\ge 0$, by the quadratic variation 
of semimartingale $\eb{t},\,t\ge 0$:  
\begin{align*}
 A_{t}:=\int _{0}^{t}e^{2B_{s}}\,ds.
\end{align*}
This exponential additive functional of Brownian 
motion appears in a number of areas of probability theory 
such as the pricing of Asian options in mathematical 
finance, and is known for its close relationship with planar  
Brownian motion (or two-dimensional Bessel process); see 
the detailed surveys \cite{mySI, mySII} by Matsumoto and Yor. 
Let $\be =\{ \be _{t}\} _{t\ge 0}$ be another one-dimensional 
standard Brownian motion which we assume to be independent of 
$B$ throughout the paper. The celebrated Bougerol's identity in law 
(\cite{bou}) then asserts that, for every fixed $t>0$, 
\begin{align}\label{;boug}
 \be _{A_{t}}\eqd \sinh B_{t}. 
\end{align}
When $t>0$, the exponential functional $A_{t}$ admits a density 
with closed-form expression 
\begin{align}\label{;expr}
 \exp \left( 
 \frac{\pi ^{2}}{8t}
 \right) 
 \ex \!\left[ \frac{\cosh B_{t}}{\sqrt{2\pi v^{3}}}
 \exp \left( 
 -\frac{\cosh ^{2}B_{t}}{2v}
 \right) \cos \left( \frac{\pi }{2t}B_{t}\right) 
 \right] ,\quad v>0; 
\end{align}
see \cite[Equation~(1.5)]{ag}. Because the expression involves 
an oscillatory integral, it sometimes is not easy to analyze; for 
instance, as $t\downarrow 0$, the prefactor in \eqref{;expr} 
grows rapidly while, in view of the Riemann--Lebesgue lemma, 
the expectation part converges to $0$, which causes difficulties 
in numerical evaluation for a small value of $t$.
Nonetheless, 
Bougerol's identity \eqref{;boug} makes the law of $A_{t}$ 
fairly tractable; as an illustration, we readily see that 
its Mellin transform is expressed as 
\begin{align*}
 \ex \!\left[ 
 \left( A_{t}\right) ^{\nu -1}
 \right] 
 =\frac{\sqrt{\pi }}{2^{\nu -1}\Gamma (\nu -1/2)}
 \ex \!\left[ |\sinh B_{t}|^{2\nu -2}\right] ,\quad \nu >1/2
\end{align*}
(see, e.g., \cite[Equation~(4.g)]{yor}), thanks to the identity.  
Here $\Gamma (\cdot )$ is the gamma function.
For a more detailed account of Bougerol's 
identity, as well as for recent progress in its study such as extensions 
to other processes, see the survey \cite{vak} by Vakeroudis. 
We also refer to \cite{har} for multidimensional extensions other than 
dealt with in the present paper.

Let $L=\{ L^{x}_{t}\} _{t\ge 0,x\in \R }$ and 
$\la =\{ \la ^{x}_{t}\} _{t\ge 0,x\in \R }$ be 
(the bicontinuous versions of) the local time processes of 
$B$ and $\be $, respectively. In \cite{bdy}, the following 
two-dimensional extension of Bougerol's identity \eqref{;boug} 
is shown: 
\begin{thm}[\cite{bdy}, Theorem~1.1]\label{;tbdy}
 For every fixed $t>0$, it holds that 
 \begin{align*}
 \left( 
 \be _{A_{t}},\,e^{-B_{t}}\la ^{0}_{A_{t}}
 \right) \eqd 
 \left( 
 e^{-B_{t}}\be _{A_{t}},\,\la ^{0}_{A_{t}}
 \right) 
 \eqd 
 \left( 
 \sinh B_{t},\,\sinh L^{0}_{t}
 \right) .
 \end{align*}
\end{thm}

We also refer the reader to \cite[Theorem 4.1]{vak}.
By virtue of the scaling property of Brownian motion, 
the first identity in law is an immediate consequence of 
the independence of $B$ and $\be $ and the identity in law 
\begin{align}\label{;jlaw}
 \left( 
 e^{B_{t}},\,A_{t}
 \right) \eqd 
 \left( 
 e^{-B_{t}},\,e^{-2B_{t}}A_{t}
 \right) ,
\end{align}
which follows from the time-reversal: 
$\{ B_{t-s}-B_{t}\} _{0\le s\le t}\eqd \{ B_{s}\} _{0\le s\le t}$.
Roughly speaking, their proof of \tref{;tbdy} is divided into 
three steps: first, observe that the validity of the claimed identity 
is reduced to showing 
\begin{align}\label{;eboug}
 \left( 
 e^{-B_{t}}|\be _{A_{t}}|,\,\la ^{0}_{A_{t}}
 \right) \eqd 
 \left( 
 \sinh |B_{t}|,\,\sinh L^{0}_{t}
 \right) ;
\end{align}
then, for each $p>0$, replace $t$ by an exponential 
random variable with parameter $p$ that is assumed to be 
independent of $B$ and $\be $; and then, with this replacement, 
show that the joint Mellin transforms  of both sides of 
\eqref{;eboug} agree for any $p>0$. While the original 
proof of \tref{;tbdy} was rather involved (in fact, in the 
course of the proof, an auxiliary standard exponential 
random variable independent of the other elements, is 
also considered to make the computation progress), 
we have noticed in \cite{mat} that the validity of the theorem 
is explained via the following two standard facts: given $t>0$, 
for every fixed $x\in \R $,
\begin{align}\label{;bougd}
 \eb{t}\sinh x+\be _{A_{t}}\eqd \sinh (x+B_{t}),  
\end{align}
and, for all $a\in \R $ and $b>0$, 
\begin{align}\label{;jloc}
 \pr \!\left( 
 B_{t}\le a,\,L^{0}_{t}\ge b
 \right) 
 =
 \begin{cases}
  \pr (b+B_{t}\le 0)+\pr (-a\le b+B_{t}\le 0) & (a\ge 0),\\
  \pr (b+B_{t}\le a) & (a<0).
 \end{cases}
\end{align}
To see how these two facts explain \tref{;tbdy}, 
for each $x\in \R $ and $y>0$, 
consider the probability 
\begin{align*}
 \pr \!\left( 
 \be _{A_{t}}\le x,\,e^{-B_{t}}\la ^{0}_{A_{t}}\ge y
 \right) .
\end{align*}  
In what follows, we denote by $\sinh ^{-1}x,\,x\in \R $, the 
inverse function of the hyperbolic sine function. 
Then, in the case $x<0$, by \eqref{;jloc} for $a<0$ and thanks 
to the independence of $B$ and $\be $, the above probability 
agrees with 
$
 \pr \!\left( 
 \eb{t}y+\be _{A_{t}}\le x
 \right) 
$, which, by \eqref{;bougd}, is equal to 
\begin{align*}
 \pr \!\left( 
 \sinh \left( 
 \sinh ^{-1}y+B_{t}
 \right) \le x
 \right) 
 &=\pr \left( 
 \sinh ^{-1}y+B_{t}\le \sinh ^{-1}x
 \right) \\
 &=\pr \!\left( 
 B_{t}\le \sinh ^{-1}x,\,L^{0}_{t}\ge \sinh ^{-1}y
 \right) \\
 &=\pr \!\left( 
 \sinh B_{t}\le x,\,\sinh L^{0}_{t}\ge y
 \right) ,
\end{align*} 
where, for the second equality, we used \eqref{;jloc} 
with $a=\sinh ^{-1}x<0$. 
The case $x\ge 0$ is treated in the same way by using \eqref{;jloc} 
for $a\ge 0$. We remark that the second identity in \tref{;tbdy} may 
also be proven by the same argument as above without appealing to 
\eqref{;jlaw}.

Identity \eqref{;bougd} was originally due to Alili and 
Gruet \cite[Proposition~4]{ag}; 
for the sake of the completeness of the paper, we give a proof 
of \eqref{;bougd} in Appendix, by means of stochastic differential 
equations (SDEs). Fact \eqref{;jloc} is a consequence of the well-known 
formula for the joint density of $B_{t}$ and $L^{0}_{t}$: 
\begin{align*}
 \pr \!\left( 
 B_{t}\in da,\,L^{0}_{t}\in db
 \right) 
 =\frac{|a|+b}{\sqrt{2\pi t^{3}}}
 \exp \left\{ 
 -\frac{(|a|+b)^{2}}{2t}
 \right\} dadb,\quad a\in \R ,\,b>0
\end{align*}
(see, e.g., \cite[Problem~6.3.4]{ks}), which may easily be 
deduced from the reflection principle of Brownian motion 
by virtue of L\'evy's theorem for Brownian local time.

In this paper, we develop further the aforementioned idea to obtain 
the following generalization of \tref{;tbdy} that extends 
\eqref{;bougd}:

\begin{thm}\label{;tmain}
For every fixed $t>0$ and $x\in \R $, we have 
\begin{equation}\label{;general}
\begin{split}
 &\left( 
 \eb{t}\sinh x+\be _{A_{t}},\,\ebm{t}\la ^{-\eb{t}\sinh x}_{A_{t}}
 \right) \\
 &\eqd 
 \left( 
 \ebm{t}\sinh x+\ebm{t}\be _{A_{t}},\,
 \la ^{-\sinh x}_{A_{t}}
 \right) \\
 &\eqd 
 \left( 
 \sinh (x+B_{t}),\,\sinh (|x|+L^{-x}_{t})-\sinh |x|
 \right) .
\end{split}
\end{equation}
\end{thm}

To our knowledge, the identity of the second coordinates, namely 
\begin{align}\label{;second}
 \la ^{\sinh x}_{A_{t}}\eqd \sinh (|x|+L^{x}_{t})-\sinh |x|
\end{align}
for every $t>0$ and $x\in \R $, seems to be new in its own right. 
As will be seen in \rref{;rsecond}, identity \eqref{;second} 
may be explained by the original Bougerol's identity \eqref{;boug} and 
by the fact that 
\begin{align}\label{;loc}
 \{ L^{x}_{t}\} _{t\ge 0}\eqd 
 \left\{ 
 \Bigl( 
 \max _{0\le s\le t}B_{s}-|x|
 \Bigr) ^{+}
 \right\} _{t\ge 0},
\end{align}
which follows from L\'evy's theorem for 
Brownian local time and the strong Markov property of 
Brownian motion. Here, for every $a\in \R $, we denote 
$a^{+}=\max \{ a,0\} $.

The rest of the paper is organized as follows: in the next section, 
preparing a lemma extending \eqref{;jloc}, we prove \tref{;tmain}, 
and in Appendix, we give a proof of \eqref{;bougd}.

\section{Proof of \tref{;tmain}}\label{;sptmain}

In the sequel, we fix $t>0$. We begin with the following lemma. 

\begin{lem}\label{;ljlocd}
For every $x\in \R $, it holds that, for all $a\in \R $ and 
$b>0$, 
\begin{equation}\label{;jlocd}
\begin{split}
 &\pr \!\left( 
 B_{t}\le a,\,L^{x}_{t}\ge b
 \right) \\
 &=
 \begin{cases}
  \pr (b+|x|+B_{t}\le 0)+\pr (x-a\le b+|x|+B_{t}\le 0) & (a\ge x),\\
  \pr (b+|x|+B_{t}\le a-x) & (a<x).
 \end{cases}
\end{split}
\end{equation}
\end{lem}

The above relation is obtained by the formula: 
\begin{align*}
 &\pr \!\left( 
 B_{t}\in da,\,L^{x}_{t}\in db
 \right) \\
 &=\frac{|a-x|+b+|x|}{\sqrt{2\pi t^{3}}}\exp \left\{ 
 -\frac{(|a-x|+b+|x|)^{2}}{2t}
 \right\} dadb,\quad a\in \R ,\,b>0
\end{align*}
(see, e.g., \cite[p.\,161, Formula~1.1.3.8]{bs}); one may also 
conclude \eqref{;jlocd} from \eqref{;jloc} by the strong Markov property 
of Brownian motion. For ease of presentation, we slightly modify 
relation \eqref{;jlocd} in such a way that, given $x\in \R $, 
for all $a\in \R $ and $b>0$, 
\begin{equation}\label{;jlocdd}
\begin{split}
 &\pr \!\left( 
 x+B_{t}\le a,\,L^{-x}_{t}\ge b
 \right) \\
 &=
 \begin{cases}
  \pr (b+|x|+B_{t}\le 0)+\pr (-a\le b+|x|+B_{t}\le 0) & (a\ge 0),\\
  \pr (b+|x|+B_{t}\le a) & (a<0),
 \end{cases}
\end{split}
\end{equation}
with which we prove \tref{;tmain}.

\begin{proof}[Proof of \tref{;tmain}]
In view of \eqref{;jlaw} and the scaling property of Brownian motion, 
it suffices to prove the identity between the first line and the third 
line of the claimed identity \eqref{;general}. To this end, we show that 
\begin{equation}\label{;target}
\begin{split}
 &\pr \!\left( 
 \eb{t}\sinh x+\be _{A_{t}}\le y,\,\ebm{t}\la ^{-\eb{t}\sinh x}_{A_{t}}\ge z
 \right) \\
 &=\pr \!\left( 
 \sinh (x+B_{t})\le y,\,\sinh (|x|+L^{-x}_{t})-\sinh |x|\ge z
 \right) 
\end{split}
\end{equation}
for any $y\in \R $ and $z>0$; once this is done, then we also 
have, for any $y\in \R $,  
\begin{align*}
&\pr \!\left( 
 \eb{t}\sinh x+\be _{A_{t}}\le y,\,\ebm{t}\la ^{-\eb{t}\sinh x}_{A_{t}}=0
 \right) \\
 &=\pr \!\left( 
 \sinh (x+B_{t})\le y,\,\sinh (|x|+L^{-x}_{t})-\sinh |x|=0
 \right) 
\end{align*}
thanks to \eqref{;bougd}, and the proof is completed. 

In the case $y<0$, the left-hand side of \eqref{;target} is 
rewritten, by the latter relation in \eqref{;jlocdd} and 
by the independence of $B$ and $\be $, as 
\begin{align}
 &\pr \!\left( 
 \eb{t}\sinh x+\be _{A_{t}}\le y,\,\la ^{-\eb{t}\sinh x}_{A_{t}}\ge \eb{t}z
 \right) \notag \\
 &=\pr \!\left( 
 \eb{t}z+\eb{t}\sinh |x|+\be _{A_{t}}\le y
 \right) \notag \\
 &=\pr \!\left( 
 \sinh \left( 
 \sinh ^{-1}(z+\sinh |x|)+B_{t}
 \right) \le y
 \right) , \label{;negl}
\end{align}
where the last line is due to \eqref{;bougd}. On the other hand, 
the right-hand side of \eqref{;target} is rewritten as 
\begin{align*}
 &\pr \!\left( 
 x+B_{t}\le \sinh ^{-1}y,\,L^{-x}_{t}\ge \sinh ^{-1}(z+\sinh |x|)-|x|
 \right) \notag \\
 &=\pr \!\left( 
 \sinh ^{-1}(z+\sinh |x|)-|x|+|x|+B_{t}\le 
 \sinh ^{-1}y
 \right) ,
\end{align*}
which agrees with \eqref{;negl}. Here we used the latter relation in 
\eqref{;jlocdd} for the equality.

The case $y\ge 0$ is treated similarly. In this case, by the 
former relation in \eqref{;jlocdd}, one sees that 
the left-hand side of \eqref{;target} is rewritten as 
\begin{align}\label{;posl}
 &\pr \!\left( 
 \eb{t}(z+\sinh |x|)+\be _{A_{t}}\le 0
 \right) 
 +\pr \!\left( 
 -y\le \eb{t}(z+\sinh |x|)+\be _{A_{t}}\le 0
 \right) \notag \\
 &=\pr \!\left( 
 \sinh \left( 
 \sinh ^{-1}(z+\sinh |x|)+B_{t}
 \right) \le 0
 \right) \notag \\
 &\quad +\pr \!\left( 
 -y\le \sinh \left( 
 \sinh ^{-1}(z+\sinh |x|)+B_{t}
 \right) \le 0
 \right) .
\end{align}
On the other hand, the right-hand side of \eqref{;target} is 
rewritten as 
\begin{align*}
 &\pr \!\left( 
 x+B_{t}\le \sinh ^{-1}y,\,L^{-x}_{t}\ge \sinh ^{-1}(z+\sinh |x|)-|x|
 \right) \\
 &=\pr \!\left( 
 \sinh ^{-1}(z+\sinh |x|)+B_{t}\le 0
 \right) \\
 &\quad +\pr \!\left( 
 -\sinh ^{-1}y\le \sinh ^{-1}(z+\sinh |x|)+B_{t}
 \le 0\right),
\end{align*}
which agrees with \eqref{;posl}. Here we used the former 
relation in \eqref{;jlocdd} for the equality. 
Therefore we have obtained \eqref{;target} for all $y\in \R $ and $z>0$ 
and the theorem is proven.
\end{proof}

\begin{rem}\label{;rsecond}
In view of \eqref{;loc}, the left-hand side of \eqref{;second} 
is identical in law with 
\begin{align*}
 \Bigl( \max _{0\le s\le A_{t}}\be _{s}-\sinh |x|\Bigr) ^{+}
\end{align*}
owing to the independence of $B$ and $\be $, which, by the 
reflection principle of Brownian motion, has the same law as 
\begin{align*}
 \left( |\be _{A_{t}}|-\sinh |x|\right) ^{+},
\end{align*}
and hence we have the identity in law 
\begin{align}\label{;secondl}
 \la ^{\sinh x}_{A_{t}}\eqd 
 \left( \sinh |B_{t}|-\sinh |x|\right) ^{+}
\end{align}
thanks to the original Bougerol's identity \eqref{;boug}. 
On the other hand, by \eqref{;loc} and the reflection principle 
of Brownian motion again, the right-hand side of \eqref{;second} 
is identical in law with 
\begin{align*}
 &\sinh \left\{ 
 |x|+\Bigl( \max _{0\le s\le t}B_{s}-|x|\Bigr) ^{+}
 \right\} -\sinh |x|\\
 &\eqd \sinh \left\{ 
 |x|+\left( |B_{t}|-|x|\right) ^{+}
 \right\} -\sinh |x|,
\end{align*}
which coincides with the right-hand side of \eqref{;secondl}.
\end{rem}

\appendix 
\section*{Appendix}
\renewcommand{\thesection}{A}
\setcounter{equation}{0}

This appendix is devoted to the proof of \eqref{;bougd}. 
We follow an inventive argument by Alili, Dufresne and Yor \cite{ady}, 
which is also reproduced in \cite[Appendix]{cmy}.

For every fixed $x\in \R $, 
define the process $X^{x}=\{ X^{x}_{t}\} _{t\ge 0}$ by 
\begin{align*}
 X^{x}_{t}:=\ebm{t}\sinh x +\ebm{t}\!\int _{0}^{t}e^{B_{s}}\,d\br _{s},
\end{align*}
with $\br =\{ \br _{t}\} _{t\ge 0}$ a one-dimensional standard 
Brownian motion independent of $B$. By It\^o's formula, 
we see that $X^{x}$ satisfies the SDE 
\begin{align*}
 dX^{x}_{t}=\sqrt{1+(X^{x}_{t})^{2}}\,d\ga ^{x}_{t}+\frac{1}{2}X^{x}_{t}\,dt,
 \quad X^{x}_{0}=\sinh x, 
\end{align*}
where we set $\ga ^{x}=\{ \ga ^{x}_{t}\} _{t\ge 0}$ by 
\begin{align*}
 \ga ^{x}_{t}:=\int _{0}^{t}\frac{-X^{x}_{s}\,dB_{s}+d\br _{s}}
 {\sqrt{1+(X^{x}_{s})^{2}}}.
\end{align*}
The process $\ga ^{x}$ is a continuous local martingale with quadratic 
variation $\langle \ga ^{x}\rangle _{t}=t,\,t\ge 0$, hence is a 
Brownian motion. Since the coefficients of the above SDE are 
Lipschitz continuous, there exists a unique strong solution, and it is 
easily checked that the solution is given by 
\begin{align}\label{;expliy}
 X^{x}_{t}=\sinh \left( x+\ga ^{x}_{t}\right) ,\quad t\ge 0. 
\end{align}
On the other hand, because of the multidimensional time-change 
theorem (see, e.g., \cite[Theorem~3.4.13]{ks}), 
we may find a one-dimensional standard Brownian motion 
$\be =\{ \be _{t}\} _{t\ge 0}$ independent of $B$, such that, a.s., 
\begin{align*}
 X^{x}_{t}=\ebm{t}\sinh x+\ebm{t}\be _{A_{t}},\quad t\ge 0. 
\end{align*}
Therefore, for every fixed $t>0$, by the scaling property of Brownian 
motion and the independence of $B$ and $\be $,
\begin{equation}\label{;ideny}
 \begin{split}
 X^{x}_{t}&\eqd \ebm{t}\sinh x+\be _{e^{-2B_{t}}\!A_{t}}\\
 &\eqd \eb{t}\sinh x+\be _{A_{t}}, 
\end{split}
\end{equation}
where the second line follows from \eqref{;jlaw}.
Comparing \eqref{;ideny} and \eqref{;expliy} entails \eqref{;bougd} 
because $\ga ^{x}$ is a standard Brownian motion.
\bigskip 

\noindent 
{\bf Acknowledgements.} 
The first author has been supported in part by JSPS KAKENHI 
Grant Number~22K03330.



\begin{thebibliography}{99}

\bibitem{ady} L.~Alili, D.~Dufresne, M.~Yor, Sur l'identit\'e de Bougerol pour les fonctionnelles exponentielles du mouvement brownien avec drift, in: Exponential Functionals and Principal Values Related to Brownian Motion: A collection of research papers, M.~Yor (ed.), pp.~3--14, Biblioteca de la Revista Matem\'atica Iberoamericana, Rev.\ Mat.\ Iberoamericana,  Madrid, 1997.

\bibitem{ag} L.~Alili, J.-C.~Gruet, An explanation of a generalized Bougerol's identity in terms of hyperbolic Brownian motion, in: Exponential Functionals and Principal Values Related to Brownian Motion: A collection of research papers, M.~Yor (ed.), pp.~15--33, Biblioteca de la Revista Matem\'atica Iberoamericana, Rev.\ Mat.\ Iberoamericana,  Madrid, 1997. 

\bibitem{bdy} J.~Bertoin, D.~Dufresne, M.~Yor, Some two-dimensional extensions of Bougerol's identity in law for the exponential functional of linear Brownian motion, Rev.\ Mat.\ Iberoam.\ {\bf 29} (2013), 1307--1324.

\bibitem{bs} A.N.~Borodin, P.~Salminen, Handbook of Brownian Motion -- Facts and Formulae, corrected reprint of 2nd ed., 2002, Birkh\"auser, Basel, 2015. 

\bibitem{bou} Ph.~Bougerol, Exemples de th\'eor\`emes locaux sur les groupes r\'esolubles, Ann.\ Inst.\ H.\ Poincar\'e Sect.\ B (N.S.) {\bf 19} (1983), 369--391.

\bibitem{cmy} A.~Comtet, C.~Monthus, M.~Yor, Exponential functionals of Brownian motion and disordered systems, J.\ Appl.\ Probab.\ {\bf 35} (1998), 255--271, also in: \cite{yorm}, pp.~182--203.

\bibitem{har} Y.~Hariya, On some identities in law involving exponential functionals of Brownian motion and Cauchy random variable, Stochastic  Process.\ Appl.\ {\bf 130} (2020), 5999--6037.

\bibitem{ks} I.~Karatzas, S.E.~Shreve, Brownian Motion and Stochastic Calculus, 2nd ed., Springer, New York, 1991.

\bibitem{mySI} H.~Matsumoto, M.~Yor, Exponential functionals of Brownian motion, I: Probability laws at fixed time, Probab.\ Surv.\ {\bf 2} (2005), 312--347. 

\bibitem{mySII} H.~Matsumoto, M.~Yor, Exponential functionals of Brownian motion, II: Some related diffusion processes, Probab.\ Surv.\ {\bf 2} (2005), 348--384. 

\bibitem{mat} Y.~Matsumura, On Spitzer's theorem for the winding number of planar Brownian motion and exponential Brownian functionals (in Japanese), Master Thesis, Tohoku University, 2022.

\bibitem{vak} S.~Vakeroudis, Bougerol's identity in law and extensions, Probab.\ Surv.\ {\bf 9} (2012), 411--437. 

\bibitem{yor} M.~Yor, On some exponential functionals of Brownian motion, Adv.\ in Appl.\ Probab.\ {\bf 24} (1992), 509--531, also in: \cite{yorm}, pp.~23--48. 

\bibitem{yorm} M.~Yor, Exponential Functionals of Brownian Motion and Related Processes, Springer, Berlin, 2001. 

\end{thebibliography}
\end{document}